\documentclass{amsart}
\usepackage{amssymb,latexsym,amsmath,amsfonts,amsthm,color}
\newtheorem{theorem}{Theorem}[section]
\newtheorem{proposition}[theorem]{Proposition}
\newtheorem{corollary}[theorem]{Corollary}
\newtheorem{definition}[theorem]{Definition}

\newtheorem{lemma}[theorem]{Lemma}
\newtheorem{example}[theorem]{Example}

\newcommand{\lk}{{\rm lk}}

\newcommand{\reals}{{\mathbb R}}

\newcommand{\R}{{\mathcal R}}

\newcommand{\id}{{\bf 1}}
\newcommand{\s}{\sigma}

\newcommand{\sm}{{\setminus}}
\newcommand{\ra}{{\rightarrow }}

\begin{document}
\title [Faces of generalized cluster complexes]
{Faces of generalized cluster complexes and noncrossing partitions}

\author{Eleni Tzanaki}
\address{Department of Mathematics\\
University of Crete\\
71409 Heraklion, Crete, Greece}
\email{etzanaki@math.uoc.gr}
\date{\today}
\thanks{The present research will be part of the author's Ph.D thesis at the University of Crete}
%
  \begin{abstract}
  Let $\Phi$ be an   finite root system with corresponding reflection group $W$
   and let  $m$ be a nonnegative integer.
    We consider the   generalized cluster complex $\Delta^m(\Phi)$  defined by S.~Fomin and N.~Reading
   and the poset $NC_{(m)}(W)$ of $m$-divisible noncrossing partitions defined by D.~Armstrong.
    We give a characterization of the faces of $\Delta^m(\Phi)$ in terms of $NC_{(m)}(W)$,
     generalizing that of T.~Brady and C.~Watt given  in the case $m=1$.
       Making use of this, we give a case free proof of  a conjecture of 
      F.~Chapoton and D.~Armstrong,
     which relates a certain refined face count of $\Delta^m(\Phi)$ with the
         M\"obius function of  $NC_{(m)}(W)$.
    \end{abstract}

\maketitle
\section{Introduction}
 Let $\Phi$ be a finite  root system of rank $n$ with associated reflection group $W$.
Let $\Phi^+$ be a positive system for $\Phi$  with corresponding simple system $\Pi$.
 Motivated by their theory of cluster algebras \cite{FZ2},  S.~Fomin and A.~Zelevinsky
 introduced the cluster complex $\Delta(\Phi)$  \cite{FZ1}. This is a pure $(n-1)$-dimensional
simplicial complex on the vertex set  $\Phi^+ \cup (-\Pi) $ which is
homeomorphic to a sphere \cite{FZ1}.
   Later, S.~Fomin and N.~Reading \cite{FR2} introduced the generalized cluster complex $\Delta^m(\Phi)$, 
 where $m$ is any nonnegative integer (see also \cite{Tz}).
 This  is a simplicial complex   on the vertex set of {\em colored almost positive roots}
 $\Phi^{m}_{\geq -1}$, that is   the set consisting of $m$ (colored) copies of each
  positive root and one copy of each negative simple root. In the case $m=1$,
  the generalized cluster complex $\Delta^m(\Phi)$  reduces to  $\Delta(\Phi)$.
 The complex $\Delta^m(\Phi)$ has remarkable properties and surprising connections with other combinatorial
  objects like  $m$-divisible noncrossing partitions $NC_{(m)}(W)$ \cite{Arm,Kr1,Kr2} and 
  Catalan hyperplane arrangements \cite{Atha2,Atha3,ATz1}.  For instance, if $\Phi$ is irreducible
   the number of facets of $\Delta^m(\Phi)$ is
    equal to the generalized Catalan number ${\rm N_m}(\Phi)=\prod_{i=1}^n \frac{e_i+mh+1}{e_i+1}$ \cite{Atha2}.
    Moreover,  the entries  of the $f$ and $h$-vector  of $\Delta^m(\Phi)$, as well as those of the natural
   subcomplex $\Delta^m_+(\Phi)$ called its {\em positive part},
   have many interesting combinatorial interpretations \cite{Atha3,ATz1}.

   In his thesis \cite{Arm} D.~Armstrong defined the poset 
   $NC_{(m)}(\gamma)$ where $\gamma$ is a Coxeter element of $W$. It is proved that the isomorphism type 
     of the poset $NC_{(m)}(\gamma)$ is independent of the Coxeter element  $\gamma$. 
    We denote this poset by $NC_{(m)}(W)$ when 
    the choice of $\gamma$ is irrelenvant and call it the poset of {\em $m$-divisible noncrossing partitions}. 
    The poset  $NC_{(m)}(W)$ is a graded meet-semilattice and  it reduces to the 
   lattice of noncrossing partitions  $NC_W$ \cite{Be1,BW1} 
    associated to $W$  in the case $m=1$. If $W$ is irreducible, 
   the number of elements of $NC_{(m)}(W)$ is equal to ${\rm N_m}(\Phi)$ 
   and  the $h$-polynomial of $\Delta^m(\Phi)$ is equal to the rank generating polynomial
  of $NC_{(m)}(W)$ \cite{Arm,FR2}.

  The first main result of this paper is a new characterization of the faces of $\Delta^m(\Phi)$
   in terms of $m$-divisible noncrossing partitions. 
   More presicely, let $\Phi$ be an irreducible root system and let $\Pi=\Pi_+ \cup \Pi_-$ be
   a partition of the set  of simple roots $\Pi$ into two disjoint sets such that the roots
      within  each are pairwise orthogonal. Let $\gamma$ be a bipartite
  Coxeter element with respect to this partition of $\Pi$ (see (\ref{bipartite})).
  Consider a  face $\s$  of $\Delta^m(\Phi)$ and let $\s^i$ be the subset of $\s$ consisting of
   positive roots of color $i$ and  $\sigma_{\pm}=\sigma \cap (-\Pi_\pm)$.
   If $\tau \subseteq \Phi^m_{\geq -1}$ such that either $\tau \subseteq (-\Pi)$ or $\tau$ consists of 
    positive roots of the same color, we denote by $w_\tau$ the product of reflections throught the roots in 
    $\tau$  taken in a certain order (see (\ref{total order})). 
   We make the convention that $w_\emptyset=\id$, where $\id$ is the  identity in $W$.
  The faces of $\Delta^m(\Phi)$ can be characterized by the following criterion.

   \begin{theorem}
  \label{main}
    The set  $\s \subseteq \Phi_{\geq -1}^m$ is a face of $\Delta^m(\Phi)$ if and only if 
     the sequence $ ( w_{\s_+} w_{\s^m},w_{\s^{m-1}},\ldots,w_{\s^2},w_{\s^1}w_{\s_-})$ 
      is an element of   $NC_{(m)}(\gamma)$ of rank $|\s|$.
    \end{theorem}

   \vspace{0.1in}
   Theorem \ref{main} specializes to the characterization of T.~Brady and C.~Watt \cite[Section 8]{BW2}
   in the case $m=1$.
   Moreover, making use of the above criterion one  may discover many interesting properties
   of $\Delta^m(\Phi)$, for instance that $\Delta^m(\Phi)$ is shellable and  $(m+1)$-Cohen-Macaulay.
    This is discussed  in the  article \cite{ATz2}.

  Our second main result is the proof of a conjecture  that relates the cluster complex
  $\Delta^m(\Phi)$ with the poset $NC_{(m)}(W)$.
   F.~Chapoton \cite{Ch} conjectured a surprising enumerative  relation between
   a refined face count  of $\Delta(\Phi)$ and the M\"obius function of $NC_W$.
  This conjecture, which was proved by Athanasiadis \cite{Atha4}, can be
  stated  in the $m \geq 1$ case as follows \cite{Arm} (see also  \cite{Kr1,Kr2}).
  The $F$-triangle for $\Delta^m(\Phi)$ is  defined by the generating function
  \begin{equation}
  \label{F-triangle}
     F^{(m)}_\Phi (x,y)= \sum_{k=0}^n \sum_{l=0}^n f_{k,l}(\Phi,m) x^k y^l,
  \end{equation}
  where $f_{k,l}(\Phi,m)$ is the number of faces of $\Delta^m(\Phi)$ consisting of
  $k$ colored  positive roots and $l$ negative simple roots.
 The $M$-triangle for $NC_{(m)}(W)$ is defined similarly as
  \begin{equation}
\label{M-triangle}
  M_W^{(m)}(x,y)= \sum_{ \tiny{ \begin{array}{c} a \leq b \\ \mbox{ in } NC_{(m)}(W) \end{array}}}
  \mu(a,b) x^{rk(b)-rk(a)} y^{rk(a)},
 \end{equation}
  where $\leq$ denotes the order relation in $NC_{(m)}(W)$,  $\mu$ stands for its M\"obius
  function  and $rk(a)$ is the rank of $a \in NC_{(m)}(W)$.

 The following relation was formulated by F.~Chapoton  \cite[Conjecture 1]{Ch} as a conjecture in the case $m=1$ 
  and was restated by D.~Armstrong for any $m \geq 1$.
\begin{theorem}
\label{conjecture}
Let $\Phi$ be a finite root system of rank $n$ with corresponding reflection group $W$ and
let $m$ be a nonnegative integer.
The $F$-triangle for $\Delta^m(\Phi)$ and the $M$-triangle for $NC_{(m)}(W)$ are related by the equality
 \begin{equation}
\label{M=F}
  (1-y)^n F^{(m)}_\Phi(\frac{x+y}{1-y},\frac{y}{1-y})=M^{(m)}_W(-x,-y/x).
\end{equation}
 \end{theorem}

 \noindent
 \vspace{0.1in}

 Theorem \ref{conjecture} has been proved in part  by C.~Kratthenthaler \cite{Kr1}, in a case by case fashion,
  for all finite root systems when $m=1$
   and for those  that do not contain a copy of $D_k$ for any $k \geq 4$, in the case $m \geq 2$.
 It was  observed by C.~Krattenthaler   that the  relation (\ref{M=F}) implies the following
   interesting reciprocity  \cite[Theorem 8]{Kr2}
  \[ y^n M^{(-m)}_W (xy,1/y) = M_W^{(m)}(x,y),\]
  for which there no intrinsic explanation yet.

  This paper is organized as follows.
After providing the necessary background,  we proceed with the proof of
 Theorem \ref{main} in Section \ref{s:criterion}. In Section \ref{s:labeling} we find an EL-labelling
    of the poset $NC_{(m)}(\gamma)$, in which the falling chains are in bijection
  with facets of $\Delta^m_+(\Phi)$. Given the  above results,
 we conclude with the proof   of Theorem \ref{conjecture} in Section   \ref{s:proof},
    generalizing that of  Athanasiadis \cite{Atha4} for the case $m=1$.

\section{Preliminaries}
\label{pre}
In this section we introduce our main objects of study and  state a few
  lemmas and theorems  necessary  in later sections.

{\bf The lattice $NC_W$ and the absolute order:}
Let $W$ be a finite Coxeter group of rank $n$ and let $T$ be the set of all reflections in $W$.
    The group $W$ is  generated by $T$  and  one can define a length function $l_T$ on $W$ 
   so  that is   $l_T(w)$ is  the smallest $k$ such that $w$ can be written as a product of 
   $k$ reflections  in $T$. We define a  partial order $\leq $ on $W$  by letting
 \[ u \leq v \mbox{ if and only if } l_T(u)+l_T(u^{-1}v)=l_T(v),\]
in other words if there exists a shortest factorization of $u$ into reflections in $T$ which is a prefix
  of such a shortest factorization of $v$.
  Note that all Coxeter elements are maximal in this partial  order of $W$.
 The group $W$ acting by conjugation gives automorphisms of this partial order, since the set of reflections
 is stable under conjugation.
   We define the  {\em noncrossing partition lattice} $NC_W$ to be an interval between the identity $\id$ and
   any Coxeter element $\gamma$. Since all Coxeter elements  are conjugate, the isomorphism type  
    of $NC_W$ is independent of $\gamma$. 
     The poset $NC_W$ is a self-dual  lattice of rank $n$ \cite{BW2}.
 If we have chosen some particular  Coxeter element $\gamma$, then  we denote
  the noncrossing partition lattice by $NC_W(\gamma)$. If there is no fear of confusion we
 just write $NC(\gamma)$ suppressing $W$ in the notation.

 For later reference we summarize a few definitions and  general facts. 

\begin{definition}{\rm \cite[Definition 3.1.1]{Arm} }
\label{minimal factors}
The $m$-tuple $(w_1,\ldots,w_m) \in W^m$ is a {\em minimal factorization} of $w \in W$ if
   $ w=w_1 w_2 \cdots w_m$ and $ l_T(w)=\sum_{i=1}^m l_T(w_i)$.
\end{definition}

\begin{lemma}
\label{conjugate-shift}
\itemsep=0pt
\item[ {\rm (i)}]
Conjugate elements have the same length.
\item[{\rm (ii)}]{\rm \cite[Lemma 3.1.2]{Arm}}
  If $(w_1,w_2,\ldots,w_k) $ is a minimal factorization of $w \in W$,
   then so is the $k$-tuple
   \[ (w_i, w_i^{-1}w_1 w_i,\ldots, w_i^{-1}w_{i-1}w_i,w_{i+1},\ldots,w_k),\]
for every $ 1 \leq i \leq k$.
 In other words, for every  $w_i$
   there is a minimal factorization of $w$ with $w_i$ in the first place.
  \item[{\rm(iii)}]
 If $(w_1,w_2,\ldots,w_k)$ is a minimal factorization of $w$,
 then $ w_{i_1} w_{i_2} \cdots w_{i_\ell} \leq w $
    and $ l_T(w_{i_1}\cdots w_{i_\ell})=\sum_{j=1}^\ell l_T(w_{i_j})$,
 for every $1\leq i_1 < i_2 <\cdots<i_\ell \leq k$. 
 \end{lemma}
  \emph{Proof. }
  Part (i) follows from the fact that all reflections in $W$ are conjugate and then part (ii) is immediate.
  To prove (iii) it suffices to  apply repeatedly (ii) for every $w_{i_j}$\; $j=1,\ldots,\ell$, 
   to get a minimal factorization of $w$ starting with $w_{i_1}w_{i_2}\cdots w_{i_\ell}$. 

  \qed

 \begin{lemma}{\rm \cite[Lemma 2.1]{Atha4}}
  \label{lemma1}
Let $a,b,w$ be elements of $W$.
\item[{\rm (i)}]$a \leq a w \leq b$ if and only if $ w \leq a^{-1} b \leq b$.
\item[{\rm (ii)}]$a \leq a w \leq b$ if and only if $ a \leq b w^{-1}  \leq b$.
\item[{\rm (iii)}]$a \leq b$ if and only if $a^{-1}b \leq b$ and, in that case,
  the interval   $[a,b]$ is isomorphic to $[\id,a^{-1}b]$.
 \end{lemma}

\begin{lemma}
\label{lemma10}
\itemsep=0pt
 \item[{\rm(i)}]{\rm \cite[Lemma 2.1(iv)]{ABW2}}
If $a,b \leq c \leq w$ for some $w \in W$ and $ab \leq w$ then $ab \leq c$.
 \item[{\rm(ii)}]{\rm \cite[Relation (3)]{BW2}}
If $a\leq b \leq c$ then $a^{-1}b \leq a^{-1}c$ and $b a^{-1} \leq c a^{-1}$.
\item[{\rm (iii)}]{\rm \cite[Relation (7)]{BW2}}
For distinct reflections $t_1$,$t_2$ and $w \in W$ we have
 $t_1 t_2 \leq w \Leftrightarrow t_2 \leq t_1 w \Leftrightarrow t_1 \leq w t_2$.
\end{lemma}

The following simple observation is immediate from the definitions and Lemma \ref{conjugate-shift}(i). 
\begin{lemma} 
\label{remark2} 
For $u,v,w \in W$ we have $u \leq v $ if and only if $w u w^{-1} \leq w v w^{-1}. $
\qed
\end{lemma}

{\bf Posets and EL-labelings:}
Let $(P,\leq)$ be a finite graded poset. We say that $y$ covers $x$ and write $x \ra y$ if
 $x<y$ and $x<z \leq y$ holds only for $z=y$.
Let $\mathcal E(P)$ be the set of covering relations in $P$ and consider $\Lambda$  a totally ordered set.
 An {\em edge labeling} of $P$ with label set $\Lambda$ is a map $\lambda:{\mathcal E}(P) \ra \Lambda$.
If $C:x_0 \ra x_1 \ra\cdots \ra x_r$ is  an unrefinable chain we let
  $\lambda(C)=(\lambda(x_0\ra x_1),\lambda(x_1\ra x_2),\ldots,\lambda(x_{r-1}\ra x_r))$ be
 the label of $C$ with respect to $\lambda$ and we call $C$ {\em rising or falling} with respect to
  $\lambda$ if the entries of $\lambda(C)$ strictly increase or weakly decrease respectively,
  in the total order of $\Lambda$.
Since $P$ is graded, all maximal chains $C$  in every interval $I=[x,y]$ have the same
 length, equal to $rk(y)-rk(x)$.
 The edge labeling allows us to order the maximal chains in $I$ by
 ordering the labels  $\lambda(C)$ lexicographically.
That is,  $C$ is lexicographically smaller than  $C'$
  if $\lambda(C)$ precedes $\lambda(C')$ in the lexicographic odrer induced by the total order of
 $\Lambda$.
 \begin{definition}{\rm (\cite{Bj})}
An edge labeling $\lambda$ of $P$ is called EL-labeling if for every non-singleton interval $[x,y]$
 in $P$ \\
  {\rm (i)} there is a unique rising maximal chain in $[x,y]$ and\\
  {\rm (ii)} this chain is lexicographically smallest among all maximal chains in $[x,y]$
   with respect to $\lambda$.
\end{definition}

%

 If $P_1,\ldots,P_m$ are posets  then their direct product $P=P_1\times \cdots \times P_m$
 is the poset on the set $\{(x_1,\ldots,x_m): x_i \in P_i\}$, such that
  $(x_1,\ldots,x_m) \leq (x'_1,\ldots,x'_m)$ in $P$ if $x_i \leq x'_i$ in $P_i$ for every $1\leq i \leq m$.
 Moreover $(x'_1,\ldots,x'_m)$ covers $(x_1,\ldots,x_m)$ in $P$  if there is some $1\leq i_0 \leq m$
  such that $x_i=x'_i$ for all $1\leq i \not= i_0 \leq m$ and $x'_{i_0}$ covers $x_{i_0}$ in $P_{i_0}$.
 The following lemma is a reformulation of \cite[Theorem 4.3]{Bj}.
 \begin{lemma}\cite{Bj}
 \label{inherit labeling}
 If $P_1,P_2,\ldots,P_m$ are graded  posets that admit an EL-labeling, then
$P=P_1\times P_2\times \cdots \times P_m$ admits an EL-labeling  as well.
 \end{lemma}
 The EL-labelling  that $P$ inherits from $P_1,\ldots,P_m$ is the following.
  Assume that $\Lambda_i$ is the totally ordered  set of the EL-labels of  $P_i$.
    If $\lambda_i \in \Lambda_i$,
  we totally order the $m$-tuples $(\emptyset,\ldots,\emptyset,\lambda_i,\emptyset,\ldots,\emptyset)$ with
  $\lambda_i$ in the $i$-th place  by letting
  \begin{equation}
  \label{cartesian-order}
 (\emptyset,\ldots,\emptyset,\lambda_i,\emptyset,\ldots,\emptyset) \leq
           (\emptyset,\ldots,\emptyset,\lambda_j,\emptyset,\ldots,\emptyset)
\mbox{ if and only if  }
   \end{equation}
  \begin{center}
  $i=j$ and $\lambda_i \leq \lambda_i $ in the total order of $\Lambda_i$ \\
  or $i > j $.
   \end{center}

   If $x=(x_1,\ldots,x_m)$ is covered by $ x'=(x'_1,\ldots,x'_m)$ in the poset $P$ 
    then   there is some $1\leq i_0 \leq m$ such that
  $x_i=x'_i$ for all $1\leq i\not= i_0 \leq m$ and $ x_{i_0} \ra x'_{i_0}$ in $P_{i_0}$.
 If $\lambda_{i_0}=\lambda(x_{i_0}\ra x'_{i_0})$ is the edge label from the EL-labeling of $P_{i_0}$
  then we label the edge $x\ra x'$ in $P$ by the $m$-tuple
  $(\emptyset,\ldots,\emptyset,\lambda_{i_0},\emptyset,\ldots,\emptyset)$
  with $\lambda_{i_0}$ in the $i_0$-th entry.
 The edge labels of $P$ we obtain this way, totally ordered as in  (\ref{cartesian-order}),
   form  an  EL-labeling of $P_1 \times P_2 \times \cdots \times P_m$.

{\bf The EL-labeling of $NC(\gamma)$:}
 Let $\Phi$ be an irreducible  root system of rank $n$ with positive part $\Phi^+$ and
 let $\Pi=\{ \s_1,\s_2,\ldots,\s_n\}$ be a choice of simple system  for $\Phi$ such that
 $ \Pi_+= \{ \s_1,\s_2,\ldots,\s_r \}$ and $\Pi_- = \{ \s_{r+1},\ldots,\s_n \}$ are orthonormal sets.
  For each $\alpha \in \Phi$ we denote by $R(\alpha)$ the reflection in $\reals^n$
 through the hyperplane orthogonal to $\alpha$.
A {\em bipartite Coxeter element} of the root system $\Phi$ is a Coxeter element $\gamma=\gamma_+\gamma_-$,
  where
  \begin{equation}
  \label{bipartite}
  \gamma_{\pm}=\prod_{\alpha \in \Pi_{\pm}}R(\alpha).
   \end{equation}
Thus, in the present case, the bipartite Coxeter element is $\gamma=R(\s_1)\cdots R(\s_n)$.

   Let  $N=nh/2$ be the number of positive roots of $\Phi^+$.
 For $1\leq i \leq 2N$ we define  $\rho_i = R(\s_1) R(\s_2)\cdots R(\s_{i-1})(\s_i)$
  where  the simple roots $\s_i$ are  indexed cyclically modulo $n$
 and we make the convention that $\rho_{-i}=\rho_{2N-i}$ for $i\geq 0$.
  The following formula, which is easily verified, appears  in \cite[Section 3]{BW2}
\begin{equation}
 \label{cases1}
 \rho_i=
 \begin{cases}
   \s_i  & \text{for  \ $i=1,\ldots,r$,} \\
  -\gamma(\s_i)& \text{for $i=r+1,\ldots,n$} \\
  \gamma(\rho_{i-n}) & \text{for  $i>n$.}
\end{cases}
  \end{equation}
Then we have
 \[\{ \rho_1,\ldots,\rho_N\}=\Phi^+, \]
 \[ \{ \rho_{N+i}: 1\leq i \leq r \}=-\Pi_+,\]
 \[ \{ \rho_{-i}: 0\leq i < n-r \}=\{ \rho_{N-i}: 0 \leq i <n-r  \}=-\Pi_-. \]

Moreover, the last $n-r$ positive roots   are a permutation of $\Pi_-$.
Let $\Phi_{\geq -1}=\Phi^+ \cup (-\Pi)$ be the set of {\em almost positive roots}.
We totally order the roots in $\Phi_{\geq -1}$ as follows:
\begin{equation}
\label{total order}
 \underbrace{ \rho_{-n+r+1} \prec \cdots \prec \rho_0}_{-\Pi_-}
 \prec \underbrace{ \rho_1 \prec \cdots \prec  \rho_N}_{\Phi^+}  \prec
  \underbrace{\rho_{N+1} \prec \cdots \prec \rho_{N+r}}_{-\Pi_+}
\end{equation}


Let $u,v \in NC(\gamma)$ such that $v$ covers $u$. Clearly $v=u t$ where $t$ a reflection in $W$.
 The {\em natural edge labeling} of the edge $u \rightarrow v$ is the reflection
  $t=u^{-1}v$.

\begin{theorem}{\rm \cite[Theorem 4.2]{ABW}}
\label{theor:shellability ABW}
If the set $T$  of positive roots $\Phi^+$ is totally ordered by {\rm (\ref{total order})}
and $\gamma= R(\s_1)R(\s_2)\cdots R(\s_n) $ then the natural edge labeling
 of $NC(\gamma)$ with label set $T$ is an EL-labeling.
\end{theorem}

{\bf The poset  $NC_{(m)}(W)$:}
 Here we recall some general facts for the poset $NC_{(m)}(W)$ of $m$-divisible noncrossing partitions. 
 We refer the reader to \cite{Arm} for more details.

\begin{definition} \cite[Definition 3.2.2]{Arm}
  \label{def:NC_m}
  Let $W$ be a reflection group  and fix a Coxeter element $\gamma$. 
  We denote by $NC_{(m)}(\gamma)$ the set of $m$-tuples $(w_1,\ldots,w_m) \in (NC(\gamma))^m$ 
    for which 
\item[ {\rm (i)}]  $w=w_1\cdots w_m \leq \gamma$
and
\item[ {\rm (ii)}] $l_T(w)=\sum_{i=1}^m l_T(w_i)$.
 \end{definition}

Recall that a subset $I$ of a poset $(P,\leq)$ is an {\em order ideal} if
 $x, y \in P$, $x\in I$ and $ y\leq x$ implies $y \in I$.
The set $NC_{(m)}(\gamma)$  is a subposet of $(NC(\gamma))^m$. Moreover it is an order ideal
   of $(NC(\gamma))^m$ with maximal elements  the minimal factorizations of $\gamma$ \cite[Lemma 3.4.3]{Arm}. 
 Thus, the set $NC_{(m)}(\gamma)$ inherits the partial order as well as the rank fuction of $(NC(\gamma))^m$. 
  More specificaly, for $(u_1,\ldots,u_m)$, $(w_1,\ldots,w_m) \in NC_{(m)}(\gamma)$ we have 
 \[ (u_1,u_2,\ldots,u_m) \leq (w_1,w_2,\ldots,w_m) \; \mbox{ if }\;  u_i \leq w_i\; \mbox{ for all  }\; 1\leq i \leq m \]
  and 
   \[ rk((w_1,\ldots,w_m)) =\sum_{i=1}^m l_T (w_i). \]

\vspace{0.1 in}
  By the fact that all Coxeter elements are conjugate and Lemma \ref{remark2} 
  we deduce that the isomorphism type of $NC_{(m)}(\gamma)$ is independent of the choice of $\gamma$. 
 Therefore, if we are not interested on the particular choice of $\gamma$ we write 
  $NC_{(m)}(W)$ and we call this the poset of  {\em $m$-divisible noncrossing partitions}. 
  The poset  $NC_{(m)}(W) $ is a ranked meet-semilattice \cite[Lemma 3.4.4]{Arm}.


{\bf The cluster complex $\Delta(\Phi)$:}
 Let $\Phi$ be an irreducible root system
   and  let $\Phi_{\geq -1} = \Phi^+ \cup (-\Pi)$  be the set of  almost positive roots.
    We define  the map $\R$ on   $\Phi_{\geq -1}$  as follows:
\begin{equation}
 \label{R}
 \R(\alpha)=
 \begin{cases}
     \gamma^{-1}(\alpha) & \text{if  $\alpha \not \in \Pi_+ \cup (-\Pi_-)$} \\
  -\alpha  &  \text{  if  $\alpha  \in \Pi_+ \cup (-\Pi_-)$ }
  \end{cases}
  \end{equation}

\begin{theorem} {\rm \cite{FZ1}}
\label{theor:compatibility}
There is a unique symmetric binary relation on $\Phi_{\geq -1}$ called {\em compatibility}
such that
\itemsep =0pt
\item[{\rm (i)}] $\alpha$ and $\beta$ are compatible if and only if $\R(\alpha)$ and $\R(\beta)$ are compatible,
\item[{\rm (ii)}] a negative simple root $-\s_i$ is compatible with a positive root $\beta$ if and only if
 the simple root expansion of $\beta$ does not involve $\s_i$.
\end{theorem}

The cluster complex $\Delta(\Phi)$ is the simplicial complex on the vertex set $\Phi_{\geq -1}$ such that
 a set $\s \subseteq \Phi_{\geq -1}$ is a face of $\Delta(\Phi)$ if every pair of roots in $\s$ is compatible.
 The maximal faces of $\Delta(\Phi)$ are called {\em clusters}.
 The simplicial complex $\Delta(\Phi)$ is pure, $(n-1)$-dimensional  and homeomorphic to a sphere \cite{FZ1}.

 {\bf The generalized  cluster complex $\Delta^m(\Phi)$:}
  The  {\em generalized cluster complex}  $\Delta^m(\Phi)$ \cite{FR2} 
  is a simplicial complex  on the vertex  set $\Phi^m_{\geq -1}$ of
  {\em  colored almost positive roots}.
   More specifically  $\Phi^m_{\geq -1} $ consists of the set $\Phi^m_{>0}$
  of $m$ (colored) copies $\alpha^1,\ldots,\alpha^m$ of each positive root $\alpha \in \Phi^+$ and
  one copy of each negative simple root. We make the convention that each negative simple root is colored by $1$.
  Thus,
   \[ \Phi^m_{\geq -1}= \Phi^m_{>0} \cup (-\Pi)^1 =
 \{ \alpha^k: \alpha \in \Phi^+, k \in \{1,\ldots,m\} \} \cup \{\alpha^1: \alpha \in -\Pi\}. \]
The faces of $\Delta^m(\Phi)$  are the subsets of $\Phi^m_{\geq -1}$ whose
  elements are pairwise  compatible in the sense we are going to describe just below.
   For every  $\alpha \in \Phi_{\geq -1}$ we define the {\em degree} $d(\alpha)$
   of $\alpha$ to be the smallest $d$ such that ${\mathcal R}^d(\alpha)$ is a negative simple root.
 Clearly, $d(\alpha)=0$ if $\alpha$ is a negative simple root.
   \begin{definition}{\rm \cite[Definition 2.1]{FR2}}
  \label{def:compatibility}
  Two colored roots $\alpha^k,\beta^l \in \Phi^m_{\geq -1} $ are called {\em $m$-compatible}
  if and only if one of the following conditions is satisfied:
  \item $\bullet$ $k>l$, $d(\alpha) \leq d(\beta) $ and the roots ${\mathcal R}(\alpha)$ and $\beta$ are compatible
   in the orininal noncolored sense (Theorem  \ref{theor:compatibility}),
 \item $\bullet$ $k<l$, $d(\alpha) \geq d(\beta) $ and the roots $\alpha$ and ${\mathcal R}(\beta)$ are compatible,
  \item $\bullet$ $k>l$, $d(\alpha) > d(\beta) $ and the roots $\alpha$ and $\beta$ are compatible,
   \item $\bullet$ $k<l $, $d(\alpha) < d(\beta) $ and the roots $\alpha$ and $\beta$ are compatible,
  and  \item $\bullet$ $k=l$  and the roots $\alpha$ and $\beta$ are compatible.
 \end{definition}

  It is immediate from the definition that $m$-compatibility  is  a symmetric relation.
  There is another equivalent way to define this relation, by introducing the $m$-analogue
  of the map $\R$, as follows.
 \begin{definition}{\rm \cite[Definition 2.3]{FR2}}
\label{def1}
 For $\alpha^k \in \Phi^{m}_{\geq -1}$,  we set

\begin{equation}
  \R_m(\alpha^k)=
\begin{cases}
   \alpha^{k+1}  & \text{if  $ \alpha \in \Phi^{m}_{>0}$  and $k < m$ } \\
      (\R(\alpha))^1  & \text{otherwise.}
\end{cases}
  \end{equation}
\end{definition}

\begin{theorem} {\rm \cite[Theorem 2.4]{FR2}}
\label{theor1}
 The compatibility relation on $\Phi^{m}_{\geq -1}$ has the following properties:
 \itemsep=0pt
 \item[{\rm (i)}] $\alpha^k$ is $m$-compatible with $\beta^l$
if and only if $\R_m(\alpha^k)$ is $m$-compatible with $\R_m(\beta^l)$,
\item[{\rm (ii)}] $(-\s_i)^1$ is $m$-compatible with $\beta^l$ if and only if
 the simple root expansion of $\beta$ does not involve $\s_i$.

Furthermore, the above conditions uniquely determine this relation.
\end{theorem}

  We denote by  $\Delta^m_+(\Phi)$ the {\em positive part} of $\Delta^m(\Phi)$,
  which is the induced subcomplex of $\Delta^m_+(\Phi)$ on the vertex set $\Phi^m_{>0}$
  of positive colored roots.
If $\Phi$ is reducible with irreducible components $\Phi_1,\ldots,\Phi_\ell$ then
  $\Phi^m_{\geq -1}= \bigcup_{i=1}^\ell (\Phi^m_i)_{\geq -1}$. We declare two roots in $\Phi^m_{\geq -1}$
compatible if they either belong to different components, or belong to the same component and are compatible within it.
Thus, the complex  $\Delta^m(\Phi)$ is the simplicial join of the complexes $\Delta^m(\Phi_i)$.

\section{Proof of Theorem \ref{main} }
 \label{s:criterion}

   Throughout this section we assume that  $\Phi$ is an irreducible root system and
 $\gamma$ a bipartite Coxeter element.
  Before proving Theorem \ref{main} we need to establish a few lemmas.

 Let $\mu:\reals^n \rightarrow \reals^n$ be the map introduced in \cite{BW2}
   with   \[ \mu(x)=2(I-\gamma)^{-1}(x)\] for all $x \in \reals^n$.
  In what follows, we denote by $\alpha \cdot \beta$ the inner product of $\alpha,\beta \in \reals^n$. 

\begin{lemma}
\label{lemma12}
\itemsep=0pt
\item[{\rm (i)}] \cite[Section 8]{BW2}
Let $\alpha, \beta \in \Phi_{\geq -1}$ with $\alpha \prec \beta$.
 The  roots $\alpha,\beta$  are compatible
 if and only if  $R(\beta) R(\alpha) \leq \gamma$.
\item[{\rm (ii)}]  {\rm \cite[Lemma 2.2]{ABW2}}
For nonparallel roots roots $\alpha,\beta$ we have $R(\alpha)R(\beta) \leq \gamma$ if and only if
 $\mu(\alpha)\cdot\beta=0$.
 \end{lemma}

\begin{lemma}{\rm \cite[Theorem 3.7]{BW2}}
 \label{lemma7}
For the total order {\rm (\ref{total order})} of roots in $\Phi_{\geq -1}$ we have 
 $\mu(\rho_j)\cdot \rho_{i-n} = - \mu( \rho_i) \cdot \rho_j$ for all $i,j$.
\end{lemma}

\begin{lemma}
 \label{lemma4}
  If  $\rho_i$ is  a positive root, then $\rho_{i-n}=\gamma^{-1}(\rho_i)$.
 \end{lemma}
 \emph{Proof. }
If $ n+1 \leq i \leq N$ the result is clear by  (\ref{cases1}).
If $1 \leq i \leq n$ then  recall  that $\rho_{i-n}=\rho_{2N+i-n}=\rho_{nh+i-n}$, 
where $N$ is the number of positive roots. 
  By the relation $\rho_i=\gamma(\rho_{i-n}) $ holding for $i>n$ and the fact that $\gamma^h=\id$ we  have
\[\rho_{n(h-1)+i}=\gamma(\rho_{n(h-2)+i})=\gamma^2(\rho_{n(h-3)+i})=\ldots
  =\gamma^{h-1}(\rho_i)=\gamma^{-1}(\rho_i), \]
 which completes our proof.
\qed

 \begin{lemma}
 \label{lemma11}
  Let $w \in W$ and $t_1,t_2,\ldots,t_k$ be distinct reflections in $W$ such that $t_i t_j \leq w$
  for every $i <j$. Then $t_1 t_2 \cdots t_k \leq w$ and $l_T(t_1t_2\cdots t_k)=k$.
  \end{lemma}
  {\em Proof. }
  We will proceed by induction, the case $k=2$ being trivial.
  We assume that the statement holds for $k$ and we  prove it for $k+1$.
     We have    $t_1\cdots t_k \leq w$ and $t_2\cdots t_{k+1} \leq w$ by the
   inductive assumption. Since $t_1 \leq t_1 \cdots t_k \leq w$ then 
    $t_2\ldots t_k \leq t_1 w$ by Lemma \ref{lemma10}(ii).
      Moreover, since  $ t_1t_{k+1} \leq w$ then  $t_{k+1} \leq t_1  w$ by part (iii) of  the same lemma.
   Setting  $a=t_2\cdots t_k$,  $b=t_{k+1}$, $c=t_1 w$ 
  and applying Lemma \ref{lemma10}(i) we get    $t_2\cdots t_{k+1} \leq t_1 w$. 
   Thus $t_1 \cdots t_{k+1} \leq w$ by Lemma \ref{lemma10}(ii). 

  To prove the statement for the length, note that since
  $t_2\cdots t_{k+1} \leq t_1 w$ then
  $l_T (t_1 w)= l_T( t_2 \cdots t_{k+1}) + l_T( t_{k+1} \cdots t_2 t_1 w)$.
   Equivalently,   $l_T(w)-1=k + l_T(w)- l_T(t_1 \cdots t_{k+1})$ by the inductive assumption and the fact that 
    $t_1\cdots t_{k+1} \leq w$. 
     Thus $l_T(t_1\cdots t_{k+1})=k+1$, which completes our proof.

 \qed

 In view of (\ref{cases1}) and the remarks following it, the total order (\ref{total order}) of the set of almost positive
  roots   $\Phi_{\geq -1}$ can be reformulated  as follows:
\begin{equation}
 \label{total order3}
 -\Pi_- \prec
\underbrace{ \Pi_+ \prec  \gamma(-\Pi_-)}_{d=1}
\prec
 \underbrace{\gamma(\Pi_+)\prec \gamma^2(-\Pi_-)}_{d=2}
\prec \gamma^2(\Pi_+)
  \prec \cdots \prec \gamma^{-1}(-\Pi_+) \prec  \Pi_- \prec -\Pi_+.
\end{equation}
 \begin{lemma}
  \label{remark}
 \itemsep=0pt
       \item[ {\rm (i)}]  If  $\alpha, \beta \in \Phi^+$ and   $d(\alpha) < d(\beta)$ then $\alpha \prec \beta$,
  \item [{\rm (ii)}] If $\alpha \in \Phi^+$, $ \beta \in \Phi^+ \setminus \Pi_+$ and $d(\beta) \leq d(\alpha)$
    then $\R(\beta) \prec \alpha$.
 \end{lemma}

  {\em Proof.}
    Part (i) is obvious from (\ref{total order3}). 
     To prove (ii), observe by (\ref{R}) and (\ref{total order3})   that if $\beta \in \Phi^+ \sm \Pi_+$ then 
     $\R(\beta)=\gamma^{-1}(\beta)$ so that $d(\R(\beta))=d(\beta)-1$. Thus 
   $d(\R(\beta)) < d(\alpha)$ and therefore $\R(\beta) \prec \alpha$ by  (i). \qed

\begin{lemma}
\label{same}
If $\alpha^\ell, \beta^k \in \Phi^m_{>0}$ are $m$-compatible then $\alpha \not = \beta$. 
\end{lemma}
  
{\em Proof.} Assume on the contrary that $\alpha^\ell, \alpha^k$ are $m$-compatible for some
  $\alpha \in \Phi^+$. By Definition \ref{def:compatibility} this implies that $\alpha,\R(\alpha)$ are 
   compatible in the noncolored sence. If $\alpha \in \Pi_+$ then  $\R(\alpha)=-\alpha$ by (\ref{R})
   and therefore  $\alpha$, $-\alpha$ are compatible which is absurd. If $\alpha \in \Phi^+ \sm \Pi_+$ 
   then  set  $\alpha = \rho_i$ for some $i \geq r+1$, so that $\R(\alpha)=\R(\rho_i)=$
    $\gamma^{-1}(\rho_i)=\rho_{i-n}$ by  (\ref{R}) and Lemma \ref{lemma4}. Moreover, 
    since $i \geq r+1$ then $\rho_{i-n} \prec \rho_i $  
      and therefore $R(\rho_i)R(\rho_{i-n}) \leq \gamma$   by Lemma \ref{lemma12}(i). 
       Hence $\mu(\rho_i)\cdot \rho_{i-n}=0$  by Lemma  \ref{lemma12}(ii) and thus
     $\mu(\rho_i)\cdot \rho_i =0 $ by Lemma \ref{lemma7}, which is  absurd since $\mu(\rho_i)\cdot \rho_i =1 $ 
    \cite[Corollary 3.3(b)]{BW2}. This  completes our proof.   \qed
  \vspace{0.1in}
  
 In what follows, for simplicity and since this does not affect the proof, 
 we omit the color 1 from the negative simple roots in $\Phi^m_{\geq -1}$.

  {\em Proof of Theorem \ref{main}:}
   We have to show that $\s \subseteq \Phi^m_{\geq -1}$ is a face of $\Delta^m(\Phi)$ if and only if 
   $w_\s= w_{\s_+}w_{\s^m} \cdots w_{\s^1}w_{\s_-} \leq \gamma$ and $l_T(w_\s)=|\s|$. 
  Recall that $\s \subseteq \Phi^m_{\geq -1}$ is a face of $\Delta^m(\Phi)$ if and only  if 
   any two roots in $\s$ are $m$-compatible. We claim that $\alpha^i,\beta^j \in \s$ are $m$-compatible 
   if and only if one of the following happens: 
    
  \begin{itemize}
 \itemsep=0pt
 \item[{\rm (i)}] $i=j$ and either $\alpha \prec \beta$ and $R(\beta)R(\alpha) \leq \gamma$,  
   or $\beta \prec \alpha$ and $R(\alpha)R(\beta) \leq \gamma$ 

 \item[{\rm (ii)}] either $\alpha \in -\Pi_+$, $\beta^j \in \Phi^m_{>0} $ 
  and $R(\alpha)R(\beta) \leq \gamma$, or $\alpha \in -\Pi_-$, $\beta^j \in \Phi^m_{>0} $ 
    and $R(\beta)R(\alpha) \leq \gamma$, 
  \item[{\rm (iii)}] $i<j$, $\alpha \not = \beta$ and  $R(\beta)R(\alpha) \leq \gamma$.
\end{itemize}
    In view of Lemma \ref{lemma11} with $w$ replaced by $\gamma$,  our claim is equivalent to the fact that 
   $w_\s \leq \gamma$ and $l_T(w_\s)=|\s|$.  

    To prove {\rm (i)}, recall by Definition \ref{def:compatibility}
 that roots of the same color are $m$-compatible if and only if they 
   are  compatible in the noncolored sense. Our claim then follows  by Lemma \ref{lemma12}(i). 
   For  {\rm (ii)}  notice that  $  \alpha \in -\Pi$ and  $ \beta^j \in  \Phi^m_{>0}$  are $m$-compatible 
   if and only  if $\alpha \in -\Pi,\beta \in \Phi^+$  are compatible in the 
   noncolored sense (Theorem \ref{theor:compatibility}(ii) and  \ref{theor1}(ii)). 
 Since roots in $-\Pi_+$ ($-\Pi_-$) succeed (precede) 
   all roots in $\Phi^+$, our claim is clear by  Lemma \ref{lemma12}(i). 

  For the proof of (iii) recall by Lemma \ref{same} that  if $\alpha^i$, $\beta^j$ 
   are $m$-compatible  then $\alpha \not = \beta$. 
   We first assume that  $d(\alpha)< d(\beta)$ so that $\alpha \prec \beta$ by Lemma \ref{remark}(i). 
   In this case, the roots  $\alpha^i, \beta^j$ are $m$-compatible if and only if  $\alpha, \beta$  are compatible. 
   By Lemma \ref{lemma12}(i) the roots $\alpha,\beta$ are compatible if and only if  
      $R(\beta)R(\alpha) \leq \gamma$, which proves our claim. 

   We next assume that $d(\alpha) \geq d(\beta)$, in which case  $\alpha^i, \beta^j$ are $m$-compatible
   if and only if  $\alpha$, $\R(\beta)$ are compatible. We distinguish cases.

   {\bf Case 1:} $\alpha \in \Phi^+$ and $\beta \in \Pi_+$. 
     By (\ref{R}) it is $\R(\beta)=-\beta \in -\Pi_+$ and since the roots in $\Phi^+$ precede those in  $-\Pi_+$ 
   we have $\alpha  \prec \ \R(\beta)$. In view of Lemma \ref{lemma12}(i) the roots 
    $\alpha, \R(\beta)$ are compatible if and only if $R(\R(\beta))R(\alpha) \leq \gamma$. 
    Since $R(\R(\beta))=R(-\beta)=R(\beta)$ we have that 
     $\alpha, \R(\beta)$ are compatible if and only if $R(\beta)R(\alpha) \leq \gamma$,  
        as desired.

  {\bf Case 2:} $\alpha \in \Phi_+$ and $\beta \in \Phi^+ \sm \Pi_+$.
    Since $d(\beta)\leq d(\alpha)$ we have that $\R(\beta) \prec \alpha$ by Lemma \ref{remark}(ii). 
   In view of Lemma \ref{lemma12}(i) the roots $\alpha,\R(\beta)$ are compatible if and only if 
     $R(\alpha)R(\R(\beta)) \leq \gamma$. 
     Set $\alpha=\rho_i$, $\beta=\rho_j$ with $i\geq 0$ and $j\geq r+1$ and 
     notice that $\R(\rho_j)=\gamma^{-1}(\rho_j)=\rho_{j-n}$ by (\ref{R}) and Lemma \ref{lemma4}. 
    We have   $R(\alpha)R(\R(\beta))\leq \gamma$ $\Leftrightarrow $ $R(\rho_i)R(\R(\rho_j)) \leq \gamma$ 
    $\Leftrightarrow $ $R(\rho_i)R(\rho_{j-n}) \leq \gamma$ $\Leftrightarrow $
    $\mu (\rho_i) \cdot \rho_{j-n} =0$ $\Leftrightarrow $ $\mu(\rho_j) \cdot \rho_i =0$ 
     $\Leftrightarrow $ $R(\rho_j)R(\rho_i) \leq \gamma$ $\Leftrightarrow $ $R(\beta)R(\alpha)\leq \gamma$    
        by Lemmas \ref{lemma12}(ii) and \ref{lemma7}. 
     Therefore $\alpha,\R(\beta)$ are compatible if and only if $R(\alpha)R(\beta) \leq \gamma$ 
      and this completes our proof.       \qed

\section{The EL-labeling and falling chains in  $NC_{(m)}(\gamma)$ }
\label{s:labeling}

In this section we describle the EL-labeling that $NC_{(m)}(\gamma)$ inherits
 from $NC(\gamma)$ and explain the relation of falling chains with respect to this EL-labeling
 to facets of $\Delta^m_+(\Phi)$. The idea for the EL-labeling was suggested by D.~Armstrong. 

  Since $NC(\gamma)$ admits an EL-labeling so does  $(NC(\gamma))^m$, by  Proposition \ref{inherit labeling}.
 Moreover, recall that  $NC_{(m)}(\gamma)$ is an order ideal in $(NC(\gamma))^m$  so that
  we can restrict the EL-labeling of the latter to the former.
  The labeling that $NC_{(m)}(\gamma)$ inherits from $(NC(\gamma))^m$ is the following.
  If  $(w_1,\ldots,w_m) \ra (w'_1,\ldots,w'_m)$ in $NC_{(m)}(\gamma)$ then there exists some
  $1\leq i_0 \leq m$ such that $w_i=w'_i$ for all $1\leq i \not = i_0 \leq m$ and
  $  w_{i_0}\leq w'_{i_0}= w_{i_0}t_{i_0}$ for some reflection  $t_{i_0}$ in $W$. 
     We label the edge $(w_1,\ldots,w_m) \ra (w'_1,\ldots,w'_m)$
    by $(\id,\ldots,\id,t_{i_0},\id,\ldots,\id)$ with $t_{i_0}$
  in the $i_0$-th entry  and we call this the {\em natural edge labeling} of $NC_{(m)}(\gamma)$.
    Let $\Lambda$ be  the set of $m$-tuples $(\id,\ldots,R(\alpha),\ldots,\id)$ where $R(\alpha)$ is a reflection
   through the root $\alpha$.  Following  (\ref{cartesian-order}), we totally order the elements of
  $\Lambda$ by letting 
  \begin{equation}
   \label{order4}
     (\id,\ldots,R(\alpha) ,\ldots,\id)  \leq    (\id,\ldots,R(\alpha') ,\ldots,\id)
  \end{equation}
    with  the reflection on the $i$-th  and $i'$-th entry respectively,
 if and only if
   \[ \begin{array}{c}
    i=i'  \mbox{ and } \alpha \preceq \alpha' \mbox{ or }\\
    i> i'. \end{array} \]

 The preceding discussion leads us to the following result. We should point out that 
 the following proposition is part of a stronger result \cite[Theorem 3.7.2]{Arm}, namely 
 that the poset $NC_{(m)}(\gamma) \cup \{ \hat 1\}$ of $m$-divisible noncrossing partitions 
 with a maximal element adjoined is EL-shellable. 
\begin{proposition}
\label{prop1}
 If $\Lambda$ is totally ordered as in {\rm (\ref{order4})} then
    the natural edge labeling of $NC_{(m)}(\gamma)$ is an EL-labeling. \qed
 \end{proposition}
\noindent

   Our next goal is to relate maximal falling chains in intervals 
 $[\hat 0,w]$  of $NC_{(m)}(\gamma)$ to facets  of certain subcomplexes of $\Delta^m_+(\Phi)$. 
  For a fixed  $w=(w_1,\ldots,w_m) \in NC_{(m)}(\gamma)$ 
   we define $\Delta^m_+(w)$ as the  induced subcomplex of $\Delta^m_+(\Phi)$ on the vertex set of
    colored positive roots  $\alpha^i \in \Phi^m_{>0}$ with $R(\alpha) \leq w_i$, $1\leq i \leq m$. 
  Since $m$-compatible roots are compatible and in view of Lemmas \ref{lemma12}(i) and \ref{lemma11}, 
   we deduce that $\Delta^m_+(w)$ is the simplicial complex whose faces are 
   the sets $\s \subseteq \Phi^m_{>0}$ with $(w_{\s^m},\ldots,w_{\s^1}) \leq (w_1,\ldots,w_m)$. 
   The complex $\Delta^m_+(w)$ is  the simplicial join of  $\Delta^1_+(w_i)$\; 
   $1\leq i \leq m$,  and therefore it is  pure of dimension  $rk(w)-1$.

  Consider  $w=(w_1,\ldots,w_m) \in NC_{(m)}(\gamma)$ with $rk(w)=k$  and write $\hat 0=(\id,\ldots,\id)$.  
    Let $C$ be  a maximal  falling chain in $[\hat 0,w]$ and observe that we can decompose 
      it into $m$ unrefinable parts  $C_i$,  such that $C_i$ is the subchain of $C$  whose 
       edge labels   have a reflection on the $i$-th entry. 
  Clearly, $C$ is the chain formed by arranging the $C_i$ one after the other. 
     Next, consider the subsets of $\Phi^+$
  \[ \s(C_i) =\{ \alpha \in \Phi^+:(\id,\ldots,R(\alpha),\ldots,\id) \mbox{ is an edge label on  } C_i \}\]
     for $1 \leq i \leq m$.
  Let  $w_{\s(C_i)}$ be  the product of reflections through the
  roots in  $\s(C_i)$ in decreasing order with respect to
  (\ref{total order}). This  is actually the order in which they appear on the edge labels of $C$.
  Thus $w_{\s(C_i)}=w_i$ and therefore $ (w_{\s(C_1)}, \ldots, w_{\s(C_m)}) \in NC_{(m)}(\gamma)$. 
   By Theorem \ref{main}, if we color the roots in each $\s(C_i)$ by $m-i+1$ we obtain 
    a face of $\Delta^m_+(\Phi)$ and in particular a facet of $\Delta^m_+(w)$. 

  Conversely, let $\s$ be  a facet of $\Delta^m_+(w)$, so that  $w_{\s^{m-i+1}}=w_i$ for all $1\leq i \leq m$. 
   Consider the set \[ \Lambda_C:=\{ (\id,\ldots,R(\alpha),\ldots,\id): \alpha \in \s_i^{m-i+1}, 
     R(\alpha) \mbox{ on the $i$-th entry}, 1\leq i \leq m \} \]  and  order its elements in decreasing order 
     with respect to (\ref{order4}). It is clear that $\Lambda_C$ is the label set of a maximal
    falling chain $C$ in $[\hat 0,w]$, obtained by reversing the procedure described in the previous
   paragraph.  We summarize the above facts in the following proposition. 
  \begin{proposition}
\label{chains}
  Let  $w=(w_1,\ldots,w_m) \in NC_{(m)}(\gamma)$ and let  $C$ be a maximal falling chain in $[\hat 0,w]$. 
  The map   which sends $C$ to $ \bigcup_{i=1}^m \s(C_i)^{m-i+1}$   is a bijection 
   from the set of maximal falling chains   of $[\hat 0,w]$ to the set of 
    facets of $\Delta^m_+(w)$. \qed 
   \end{proposition}  

\begin{corollary}
 \label{particular}
The set of  maximal falling chains in $NC_{(m)}(\gamma)$ bijects to  the set of facets of $\Delta^m_+(\Phi)$.
\end{corollary}
{\em Proof.} The set of maximal falling chains in $NC_{(m)}(\gamma)$ is the disjoint union of the sets of 
   such chains within each $[\hat 0,w]$ with $rk(w)=n$. Moreover, the set of facets of $\Delta^m_+(\Phi)$ 
   is the disjoint union of the set of facets of each $\Delta^m_+(w)$ where $rk(w)=n$. 
  Our claim then follows from Proposition \ref{chains}. \qed 
\vspace{0.1 in}

The reader is invited to verify Corollary \ref{particular} in the following example.
\vspace{-0.1in}
\begin{example} {\rm
 Consider the root system $A_2$ with positive  roots $\s_1$, $\s_2$ and $ \alpha=\s_1+\s_2$,  
  where $\s_1,\s_2$ are simple roots. 
 We set 
 $\Pi_+=\{\s_1 \}$, $\Pi_-=\{\s_2\}$  so that  $\s_1 \prec \alpha \prec \s_2$.  
We represent the elements of the reflection group 
  $W_{A_2}$ as permutations in $S_3$.
 The  reflection through $\s_1,\s_2$ and $ \alpha$  acts on $\mathbb R^3$ by transposing
  coordinates, so that $R(\alpha_1)=(12)$, $R(\alpha_2)=(23)$ and $R( \alpha)=(13)$.
   The bipartite Coxeter element in this case is
  $\gamma=(123)$.
 Figure {\rm \ref{figure1}} illustrates the natural edge labeling of $NC_{(2)}(\gamma)$.
  The set of maximal  falling chains in $NC_{(2)}(\gamma)$ bijects
    to the set of  positive clusters as shown in the following table.
\begin{center}
\begin{tabular}{| l| l |} \hline
 $ \mbox{\rm Edge labels of falling chains } $ &  $ \mbox{\rm positive  clusters }$ \\ \hline \hline
$(\id,(23))\ra(\id,(13))$ & $\s_2^1,\alpha^1$ \\ \hline

$(\id,(13))\ra (\id,(12))$ & $\alpha^1,\s_1^1$ \\ \hline

$((23),\id)\ra(\id,(13))$ & $\s_2^2,\alpha^1$ \\ \hline

$((23),\id)\ra ((13),\id)$ & $\s_2^2,\alpha^2$ \\ \hline

$((12),\id)\ra (\id,(23))$ & $\s_1^2,\s_2^1$ \\ \hline

$((13),\id)\ra ((12),\id)$ & $\alpha^2,\s_1^2$ \\ \hline

$((13),\id)\ra (\id,(12))$ & $\alpha^2,\s_1^1$ \\ \hline
\end{tabular}
\end{center} }
\end{example}

\vspace{0.2 in}

  \begin{lemma}
  \label{lemma14}
  For $w=(w_1,\ldots,w_m) \in NC_{(m)}(\gamma)$ the number $(-1)^{rk(w)} \mu(\hat0,w)$ is
  equal to the number of  facets  of $\Delta^m_+(w)$, 
\end{lemma}
 {\em Proof.} A standard fact on M\"obius functions of EL-shellable 
  posets is that  the M\"obius function on every  interval $[a,b]$  is equal to $(-1)^{rk(b)-rk(a)}$ 
   times the number of maximal falling chains 
  in $[a,b]$ \cite[Section 3.13]{St}. Applying this in our case with $a=\hat 0$ and $b=w$, we deduce  that 
   $(-1)^{rk(w)} \mu(\hat0,w)$ is equal to the number 
  of maximal falling chains in $[\hat 0,w]$  which, in view of Proposition \ref{chains},  
   is equal to the number of facets of $\Delta^m_+(w)$.    \qed
 
\section{Proof of Theorem \ref{conjecture}}
\label{s:proof}

 As pointed out in \cite[Proposition F(i)]{Kr1}, for any finite root systems $\Phi$ and $\Phi'$
  we have \[ F^{(m)}_{\Phi \times \Phi'}(x,y) = F^{(m)}_\Phi (x,y) \cdot F^{(m)}_{\Phi'}(x,y). \]
  Moreover, by the multiplicativity  of the M\"obius function, the above relation
  holds for the $M$-triangle as well. Therefore, it suffices to prove Theorem \ref{conjecture}
  in the case where $\Phi$ is irreducible.

The $h$-polynomial of an abstract $(n-1)$-dimensional simplicial complex $\Delta$ is defined as
 \[ h(\Delta,y)=\sum_{i=0}^n f_i(\Delta)y^i(1-y)^{n-i}, \]
 where $f_i(\Delta)$ is the number of faces of $\Delta$ of dimension $i-1$.
The link of a face $\s$ of $\Delta$ is the abstract simplicial complex
  $\lk_\Delta (\s)=\{ \tau \setminus \s: \s \subseteq \tau \in \Delta\}. $
One can check by comparing \cite[Figure 3.4]{Arm} and \cite[Theorem 9.2]{FR2}
 that for any irreducible root system $\Phi$ the $h$-polynomial of $\Delta^m(\Phi)$ coincides with
  the rank generating polynomial of $NC^{(m)}(\gamma)$, thus
 \begin{equation}
  \nonumber
   h(\Delta^m(\Phi),y)= \sum_{w\in NC_{(m)}(\gamma)} y^{rk(w)}.
  \end{equation}

The following lemma  will be used as in \cite{Atha4}.
 For $\alpha \in \Pi$ we denote by $\Phi_{\alpha} $ the standard parabolic
  root subsystem obtained by intersecting $\Phi$ with the linear span of $\Pi \setminus \{\alpha\}$.

\begin{lemma}
\label{lemma15}
Let $\Phi$ be irreducible, $\alpha \in \Pi$ and $\s \subseteq \Phi^m_{\geq -1}$.
\item[{\rm (i)}] For $\s \in \Delta^m(\Phi)$ we have $-\alpha \in \s$ if and only if
$\s \setminus \{-\alpha\} \in \Delta^m(\Phi_\alpha)$.
\item[{\rm (ii)}] For any $\beta^\ell \in \Phi^m_{>0} $ there exists $i$ such that $\R_m^i(\beta^\ell)\in (-\Pi)$.
\item[{\rm (iii)}] $\s \in \Delta^m(\Phi)$ if and only if $\R_m(\s) \in \Delta^m(\Phi)$.
\end{lemma}
{\em Proof.}
Part (i) can be verified from  \cite[Theorem 2.7]{FR2}. 
  Part (ii) follows from Definition \ref{def1} and the fact that for every $\beta \in \Phi^+$ 
  there exists some $j$ such that $ \R^j(\beta) \in (-\Pi)$ \cite[Section 3]{FZ1}. 
 Finally, part (iii) is clear from Theorem \ref{theor1}(i). \qed 

\vspace{0.1 in}
We continue with some technical lemmas required in the proof of  Theorem \ref{conjecture}.
 \begin{lemma}
  \label{link-lemma} 
  If  $\s$ is  a face of $\Delta^m(\Phi)$ and 
  $w_\s = w_{\s_+} w_{\s^m} w_{\s^{m-1}}\cdots w_{\s^1} w_{\s_-}$ then
  \[ h(\lk_\Delta(\s),y)= \sum_{ \alpha \in NC_{(m)}(\gamma w_\s^{-1})} y^{rk(\alpha)}. \]
  \end{lemma}
 \noindent
 \emph{Proof.}
   The proof is analogous to the proof of \cite[Lemma 2.6]{Atha4}, replacing $\R$ by $\R_m$
   and using Lemma \ref{lemma15}.
  \qed

\begin{lemma}
\label{interval}
If  $ (a_1,\ldots,a_m) \leq (b_1,\ldots,b_m)$ in $ NC_{(m)}(\gamma)$ then 
  $(a_1^{-1}b_1,\ldots,a_m^{-1}b_m) \in NC_{(m)}(\gamma)$ and 
  the intervals  $I=[(a_1,\ldots,a_m),$ $(b_1,\ldots,b_m)]$ and 
  $I'=[(\id,\ldots,\id),$ $(a_1^{-1}b_1,\ldots,a_m^{-1}b_m)]$ are isomorphic.  
\end{lemma}
{\em Proof.} 
It follows from Lemma \ref{lemma1}(iii)  that  since   $a_i \leq b_i$ then  $  a_i^{-1}b_i \leq b_i$. 
   Hence $(a_1,a_1^{-1}b_1,\ldots,a_m, a_m^{-1}b_m)$ is a minimal factorization of $b_1 \ldots b_m$ 
   and in view of Lemma \ref{conjugate-shift}(iii) we have that  $(a_1^{-1}b_1,\ldots,a_m^{-1}b_m)\in NC_{(m)}(\gamma)$. 
  This proves the first part of our statement. 
   
     Next, we claim  that the map   $\varphi:I \rightarrow  I'$ with 
     $\varphi((w_1,\ldots,w_m))=$  $(a_1^{-1}w_1,\ldots,$ $a_m^{-1}w_m)$ is an order preserving bijection 
       and thus the intervals $I$ and $I'$ are isomorphic. 
   We first have to check that $\varphi$ is well defined. Indeed, if $(w_1,\ldots,w_m) \in I$ then 
 $(a_1,a_1^{-1}w_1,\ldots,a_m,a_m^{-1}w_m)$ is a minimal factorization of $w_1\cdots w_m$ 
   and therefore  $(a_1^{-1}w_1,\ldots,a_m^{-1}w_m) \in NC_{(m)}(\gamma)$ by Lemma \ref{conjugate-shift}(iii). 
 Conversely,   let  $(w_1,\ldots,w_m) \in I'$ so that $w_i \leq a_i^{-1} b_i \leq b_i$ for $1\leq i \leq m$. 
In view of Lemma \ref{lemma1}(i) we have $a_i \leq a_i w_i \leq b_i$ and therefore 
    $(a_1 w_1, w_1^{-1} a_1^{-1}b_1,$ $ \ldots,$ $ a_m w_m, w_m^{-1}a_m^{-1}b_m)$ is a minimal 
  factorization of $b_1\cdots b_m$. Thus $(a_1w_1,\ldots,a_m w_m) \in NC_{(m)}(\gamma)$
     as in the previous situation.  
  That $\varphi$ is a bijection is immediate from the fact that 
   $a_i \leq w_i \leq b_i$ if and only if
   $ a_i^{-1} w_i \leq a_i^{-1}b_i$ (Lemma \ref{lemma1}(i)). Moreover,  $ rk(\varphi(w))= rk(w)-rk(a)$ 
  and therefore $\varphi$ is order preserving. This completes our proof. \qed 
\vspace{0.1 in}

\begin{lemma}
\label{lemma13}
Let  $(w_1,\ldots,w_m) \in$ $  NC_{(m)}(\gamma)$. 
There is a rank preserving bijection between the sets 
$ A=\{ (a_1,\ldots,a_m): (a_1,\ldots,a_m) \leq (a_1w_1,\ldots,a_m w_m)$ $ \mbox{ in } NC_{(m)}(\gamma) \}$ 
  and $NC_{(m)}(\gamma w_m^{-1}\cdots w_1^{-1})$.  
\end{lemma}
{\em Proof.}
    Let $a'_i= w_1 \cdots w_{i-1} a_i w_{i-1}^{-1} \cdots w_1^{-1}$ for $1\leq i \leq m$ 
  and consider the map $\phi(a_1,\ldots,a_m)=(a'_1,\ldots,a'_m)$. We will prove that $\phi$ is a rank 
  preserving bijection between the sets $A$ and $NC_{(m)}(\gamma w_m^{-1}\cdots w_1^{-1})$.
   Clearly, the map  $\phi$ as well as its inverse  $\phi^{-1}(a'_1,\ldots,a'_m)=(a_1,\ldots,a_m)$ 
   with $a_i=w_{i-1}^{-1} \cdots w_1^{-1} a'_i w_1 \cdots w_{i-1}$ are injective.  Moreover, 
  since  $a_i, a'_i$ are conjugate  then $l_T(a_i)=l_T(a'_i)$ and therefore $\phi$ is rank preserving. 
 So, it suffices to prove that $\phi$  is well defined.  
   To this end let $(a_1,\ldots,a_m) \in A$ and note that 
  $a'_1\cdots a'_m=(a_1w_1\cdots a_mw_m)w_m^{-1}\cdots w_1^{-1}$. 
  Since $a_i \leq a_i w_i$ then  $l_T(a_iw_i)=l_T(a_i)+l_T(w_i)$ and therefore 
    $(a_1,w_1,\ldots,a_m,w_m)$ is a minimal factorization  of $a_1w_1 \cdots a_mw_m$. 
 Thus, by Lemma \ref{conjugate-shift}(iii) we have
   \[ w_1 \cdots w_m \leq a_1 w_1 \cdots a_m w_m \leq \gamma \] 
  with $l_T(w_1\cdots w_m)=\sum_{i=1}^m l_T(w_i)$. In view of Lemma \ref{lemma10}(ii)
  this implies that 
 \[ (a_1w_1 \cdots a_m w_m) w_m^{-1}\cdots w_1^{-1} \leq \gamma w_m^{-1}\ldots w_1^{-1}\] 
  with 
$l_T((a_1w_1 \cdots a_m w_m) w_m^{-1}\cdots w_1^{-1})=\sum_{i=1}^m l_T(a_i)$. Equivalently, 
  $a'_1 \cdots a'_m  \leq \gamma w_m^{-1}\cdots w_1^{-1}$ and 
  $l_T(a'_1\cdots a'_m)= \sum_{i=1}^m l_T(a_i)=\sum_{i=1}^m l_T(a'_i)$, which shows 
  that indeed $(a'_1,\ldots,a'_m) \in NC_{(m)}(\gamma w_m^{-1}\cdots w_1^{-1})$.

For the converse let $(a'_1,\ldots,a'_m) \in NC_{(m)}(\gamma w_m^{-1}\cdots w_1^{-1})$
 and notice that \\ $ a_1 w_1 \cdots $ $a_m w_m= $ $ a'_1\cdots a'_m w_1 \cdots w_m$. 
 Since $a'_1 \cdots a'_m \leq \gamma w_m^{-1} \cdots w_1^{-1} \leq \gamma$ then 
  by Lemma \ref{lemma1}(ii) we have $a'_1\cdots a'_m w_1 \cdots w_m \leq \gamma$ with 
 $l_T(a'_1\cdots a'_m w_1 \cdots w_m)= \sum_{i=1}^m l_T(a'_i)+\sum_{i=1}^m l_T(w_i)$. 
 This forces $(a_1w_1,\ldots,a_m w_m) \in NC_{(m)}(\gamma)$.  
  We next prove that $(a_1,\ldots,a_m) \leq (a_1w_1,\ldots, a_mw_m)$. 
  We have $ a'_i \leq \gamma w_m^{-1}\cdots w_1^{-1}$ or equivalently 
  $ w_1 \cdots w_{i-1} a_i w_{i-1}^{-1} \cdots w_1^{-1} \leq  \gamma w_m^{-1}\cdots w_1^{-1}$. 
  In view of Lemma \ref{remark2} this implies that 
   $ a_i \leq w_{i-1}^{-1}\cdots w_1^{-1} \gamma w_m^{-1}\cdots w_i^{-1}.$
  By Lemma \ref{remark2} and elementary calculations one may check that 
   $ w_{i-1}^{-1}\cdots w_1^{-1} \gamma w_m^{-1}\cdots w_i^{-1}$  $\leq w_{i-1}^{-1} \cdots w_1^{-1} \gamma w_i^{-1}$ 
  $\leq \gamma w_i^{-1}$. Thus $a_i \leq \gamma w_i^{-1} \leq \gamma$ 
   and therefore $a_i \leq a_i w_i \leq \gamma$ by Lemma \ref{lemma1}(ii). 
  This completes our proof.  \qed

\vspace{0.1 in}
 {\em Proof of Theorem {\rm \ref{conjecture}}.}
  To simplify notation let us write $\Delta^m$, $\Delta^m_+$, $F(x,y)$ and $M(x,y)$  instead of
   $\Delta^m(\Phi)$, $\Delta^m_+(\Phi)$, $F^{(m)}_\Phi(x,y)$ and $M^{(m)}_W(x,y)$ respectively.
   We use the relation  
\begin{equation}
\label{eq3}
 (1-y)^n F(\frac{x+y}{1-y},\frac{y}{1-y})=  \sum_{\s \in \Delta^m_+}x^{|\s|} h(\lk_{\Delta^m}(\s),y),
\end{equation}
    which appears  in the course of the proof of \cite[Theorem 1.1]{Atha4} and can be generalized 
   straightforward  in the case where $m\geq 1$. 
     We have
  \[ M(-x,-y/x)= \sum_{{\tiny \begin{array}{c} a \leq b \\ a,b \in NC_{(m)}(\gamma) \end{array} }}
    \mu(a,b)(-x)^{rk(b)-rk(a)} y^{rk(a)}.\]
  Let  $ \hat 0 =(\id,\ldots,\id)$, $w=(w_1,\ldots,w_m)=(a_1^{-1}b_1,\ldots,a_m^{-1}b_m)$
  and note that $rk(w)=rk(b)-rk(a)$.
   In view of Lemmas  \ref{interval} and  \ref{lemma13}  the last sum becomes
 \[ \sum_{w \in NC_{(m)}(\gamma)} \mu(\hat 0,w)(-x)^{rk(w)}
             \sum_{a' \in NC_{(m)}(\gamma w_m^{-1}\cdots w_1^{-1})} y^{rk(a')},  \]
  which, by Lemma \ref{lemma14} is equal to 
  \[ \sum_{w \in NC_{(m)}(\gamma)} \sum_{\tiny \begin{array}{c} \s \in \Delta^m_+(w)\\ |\s|=rk(w)\end{array}}x^{rk(w)}
             \sum_{a' \in NC_{(m)}(\gamma w_m^{-1}\cdots w_1^{-1})} y^{rk(a')}. \]
  If $\s$ is a facet of $\Delta^m_+(w)$ then $w_\s=w_{\s^m}\cdots w_{\s^1}=w_1\cdots w_m$ and thus 
  the last expression is equal to 
 \[ \sum_{w \in NC_{(m)}(\gamma)} \sum_{\tiny \begin{array}{c} \s \in \Delta^m_+(w)\\ |\s|=rk(w)\end{array}}x^{|\s|}
             \sum_{a' \in NC_{(m)}(\gamma w_\s^{-1})}y^{rk(a')}, \]
  or, in view of Lemma \ref{link-lemma},  
 \[ \sum_{w \in NC_{(m)}(\gamma)} \sum_{\tiny \begin{array}{c} \s \in \Delta^m_+(w)\\ |\s|=rk(w)\end{array}}x^{|\s|}
            h(\lk_{\Delta^m}(\s),y).\]
   Observe that the sets of facets of the subcomplexes $\Delta^m_+(w)$ for $w\in NC_{(m)}(\gamma)$ form a partition 
 of $\Delta^m_+$, so that  the last sum is  
   \[ \sum_{ \s \in \Delta^m_+}  x^{|\s|} h(\lk_{\Delta^m}(\s),y)\]
     which, in view of (\ref{eq3}), completes our proof.

 \qed

{\em Aknowledgments:} I would like to thank Drew Armstrong for making available his thesis \cite{Arm}
   at an early stage  and Christos Athanasiadis for usefull comments and discussions.

\end{document}